\documentclass[11pt]{amsart}

\usepackage{latexsym}
\usepackage{amsfonts}
\usepackage{graphicx}
\usepackage{amssymb}
\usepackage{amsmath}
\usepackage{amssymb}

\newtheorem{example}{Example}

\newtheorem{definition}{Definition}

\def\classification{\@ifnextchar [{\@xfootnotenext}%
{\begingroup\let\protect\noexpand \xdef\@thefnmark{}\endgroup
\@footnotetext}}

\begin{document}


\title{Combinatorial group theory and public key
cryptography}

\author[V.~Shpilrain]{Vladimir Shpilrain}

\author[G.~Zapata]{Gabriel Zapata}

\medskip

\begin{abstract}
\noindent After some excitement generated by recently suggested
public key exchange protocols due to Anshel-Anshel-Goldfeld and
Ko-Lee et al., it is a prevalent opinion  now that the conjugacy search
problem is unlikely to provide sufficient level of security 
 if a braid group is used as the
platform. In this paper we   address the following questions: (1)
whether choosing a different group, or a class of groups, can
remedy the situation; (2) whether some other ``hard" problem from
combinatorial group theory can be used, instead of the conjugacy
search problem, in a public key exchange protocol. Another
question that we address here, although somewhat vague, is likely
to become a focus of the future research in public key
cryptography based on symbolic computation: (3) whether one can
efficiently disguise an element of a given group (or a semigroup)
by using defining relations.
\end{abstract}

\maketitle

\tableofcontents

\section{Introduction}

One of the possible generalizations of the {\it discrete logarithm
problem} to arbitrary groups is the so-called
 {\it  conjugacy search  problem}: given two elements $a, b$ of a group $G$
 and the information that $a^x=b$
for some $x \in G$, find at least one particular element $x$ like
that. Here $a^x$ stands for $xax^{-1}$. The (alleged)
computational difficulty of this  problem in some particular
groups (namely, in braid groups) has been  used in several group
based cryptosystems, most notably in \cite{AAG}  and
\cite{KLCHKP}. However, after some initial excitement (which has
even resulted in naming a new area of ``braid group cryptography"
--- see \cite{Lipmaa}, \cite{Dehornoy}),
 it seems now  that
the conjugacy search problem in a braid group cannot provide 
sufficient level of security; see \cite{Shpilrain} for
explanations.

 Therefore, one faces the following two natural questions:

\medskip

\noindent  {\bf Question 1.} Is there   a  group, or a class of
groups, where the public key exchange protocol suggested in
\cite{AAG} would be secure enough to be used in real-life
applications?

\medskip

\noindent  {\bf Question 2.} Is there another ``hard" problem in
combinatorial group theory that can be used, instead of the
conjugacy search problem, in a public key exchange protocol?
\medskip

 Without a positive answer to at least one of these questions,
it is unlikely that combinatorial group theory  will have a
significant impact on public key cryptography, which is now
dominated by methods and ideas from number theory.

 We point out one more question, which has not been getting
sufficient attention so far, but is likely to become a focus of
the future research in public key cryptography based on symbolic
computation:

\medskip

\noindent  {\bf Question 3.} Can one  efficiently disguise an
element of a given group (or a semigroup) by using defining
relations?
\medskip

 Disguising an element before transmission is sometimes called
``diffusion" --- see e.g. \cite{Garrett}. The importance of this
is rather obvious: if, for example, one transmits a conjugate
$xax^{-1}$ of a public element $a$ ``as is", i.e., without
diffusion, then the opponent can determine the private element $x$
just by inspection. Similar problem arises in any other public key
exchange protocol. In  protocols based on ideas from number
theory, the diffusion is usually provided  ``automatically", due
to various properties of the decimal or other numerical system
that is used. For instance, in the product 7$\cdot$ 3=21, the
factors 7 and 3 cannot be determined just by inspection; this is
provided simply by the way we multiply integers in the decimal
system, or, equivalently, by the existence of a simple ``normal 
form" for integers. 

 In abstract groups, we usually do not have this facility. In fact,
in an abstract group (or a semigroup), the result of
multiplication is simply concatenation:  $a\cdot b =ab$, i.e., an
extra effort is  always required to disguise factors in a
product. This is why a diffusion mechanism is of paramount
importance in any public key exchange protocol based on symbolic
computation.

 We note here that  recent work of
Myasnikov and  Ushakov \cite{Ushakov} makes it appear likely that,
speaking somewhat informally, in a ``generic" group, the amount of
work needed to disguise a ``generic" element by using defining relations is
about the same as needed to recover an element from its disguised
form. This, of course, is unacceptable in cryptographic
applications.
 It seems that the difficulty in
 disguising an element of a group (or a semigroup) by using defining relations
might be a major obstacle for using symbolic computation in public
key cryptography, and the problem of diffusion will therefore take
the center stage in   future research.

 In this paper, we contribute toward a solution of this problem in
Section \ref{diffusion} by breaking down defining relations of a
group into ``small pieces". More formally, we replace a given
group by an isomorphic group where all relators have length at
most 3. Intuitively, diffusion should be easier to achieve in
groups with shorter defining relations, so we hope that our idea
can be useful.

 As far as other questions are concerned,
 we have to say up front that, in our opinion, Question 1 has a
smaller chance for a positive answer, and it is unlikely that the
 conjugacy search problem  will be used in real-life implementations.
 Nevertheless, we study Question 1 here in Section \ref{small cancellation}
  by exploring the idea
of using random groups from  a sufficiently large {\it  class of
groups} instead of a single group. Technically, braid groups, too,
are a class of groups, but this class is too narrow in the sense
that, informally speaking,
 any (meaningful)  algorithm that works for a
particular group $B_n, n \ge 5$, would also  work for $B_m$ for
any $m \ge 5$. We may call such a class of groups
``algorithmically homogeneous". Here we draw attention to a more
diverse class of {\it  small cancellation groups} that satisfy
small cancellation conditions $C(4)$, $T(4)$, but not
$C'(\frac{1}{6})$ (see \cite{L-S}). The latter is needed to try to
avoid hyperbolic groups (all finitely presented $C'(\frac{1}{6})$
groups are hyperbolic), where the conjugacy search problem can be
solved very quickly (see \cite{KMSS1}  and  \cite{KMSS2} for
discussion).

 In the class of groups with small cancellation conditions
$C(4)$ and $T(4)$, the word problem is solvable in quadratic time
(see \cite[Theorem V.6.3]{L-S}), which meets the necessary condition
for an efficient common key extraction by authorized parties. We
note in passing that the   existence of a unique {\it normal form}
for elements of a particular group $G$ is not necessary for common
key extraction, as observed in  \cite{AAG}. If Alice and Bob have
arrived at a point where Alice has an element, say, $u$,  and Bob
has an element $v$ such that $u=v$ in $G$, then they can establish
a common key as follows. Alice chooses, privately, a finite binary
sequence $b_1, b_2, ...$, which is going to be her common secret
key with Bob. She then transmits a sequence of group elements
$u_1, u_2, ...$ such that $u_i=u$ in $G$ if  and only if $b_i=1$.
Bob  recovers the sequence $b_1, b_2, ...$ by comparing $u_1, u_2,
...$ to his $v$.

 We note that there is no known polynomial time algorithm for
solving the conjugacy search problem in an arbitrary group with
small cancellation conditions $C(4)$ and $T(4)$.

 In Section \ref{sec Alg PKC}, we consider a different problem from
combinatorial group theory that can be used in a public key
exchange protocol. This is yet another generalization of the
discrete logarithm problem.  Given a group $G$ with the semigroup
of endomorphisms $End ~G$, suppose there are two subsemigroups, $A
\subseteq End ~G$ and  $B \subseteq End ~G$, such that for any
$\alpha \in A$  and $\beta \in B$, one has $\alpha \beta = \beta
\alpha$. Let $w \in G$ be a  public element.  Then the key
exchange protocol is quite  standard: Alice chooses, privately,
some $\alpha \in A$ and sends $\alpha(w)$ to Bob. Bob chooses some
$\beta \in B$ and sends $\beta(w)$ to Alice. Since $\alpha \beta =
\beta \alpha$, both end up with a common private key
$\alpha(\beta(w))=\beta(\alpha(w))$.

 The point is, of course, in selecting a platform group $G$ and
semigroups $A, B \subseteq End ~G$ wisely, so that the
corresponding key exchange protocol is both secure and efficient.
One special case of such arrangement appears in \cite{KLCHKP},
where $G$ is a braid group $B_n$, and $A, B$ consist of  inner
automorphisms (i.e., conjugations).
  This arrangement however makes the cryptosystem vulnerable to so-called
``length based" attacks (see e.g. \cite{GKTTV}, \cite{HS},
\cite{HT}) because applying a generic automorphism to a generic
element of a group tends to increase the length of (the normal
form of) this element. To avoid attacks of this kind, we suggest
here using non-injective endomorphisms; the effect of such an
endomorphism on the length of an element is no longer predictable.

  Again,  in Section \ref{classes},
we suggest using a large class of groups instead of a fixed group
and selecting a random group from this class every time one wants
to initiate a public key exchange protocol. A particular class of
groups that we consider here is the class of {\it Artin groups of
extra large type}. Groups in this class are known to be automatic
\cite{Peifer}, which implies, in particular, that the word problem
in any group from this class is solvable in quadratic time.
Further details are given in Section \ref{Class of Artin Groups}.

Finally, we note that, as a further generalization, one can use
arbitrary well-defined mappings $\alpha, \beta$ (not necessarily
endomorphisms) of a group $G$ in the above context. A simple
example of that kind was given in \cite{Ko}; see also our Section
\ref{sec Alg PKC}.


\section{Algebraic public-key cryptographic systems}
\label{sec Alg PKC}

  The central requirement for an operational public-key
cryptographic system (PKC) is a \emph{one-way function}; in
theory, it is the security core in the development and
implementation of public-key cryptographic protocols. Let  $S$
and $T$ be two sets. In essence, a one-way function is
  a   feasibly computable function $f: S 
\rightarrow T$ such that given the \emph{image} $y=f(x)$, it is
computationally infeasible to determine a \emph{preimage} $x \in
S$.

 For an algebraic characterization of a one-way function, we
assume $S$ and $T$ to be   associative
algebraic structures   with a single binary operation,
e.g., semigroups. We  call these structures \emph{platforms} when
used in the context of  cryptography.

 Let the pair $\langle X ; R \rangle$ be a \emph{presentation} of a
semigroup $S$, where $X=\{x_1,x_2,\ldots\}$ is   a set of
generators of $S$ and $R=\{r_1=r_1',r_2=r_2',\ldots\}$   a set
of defining relations. The \emph{full transformation
semigroup} of $S$, denoted by $\mathcal{T}_S$, is the set of all
functions $S \rightarrow S$ closed under composition, see e.g.
\cite{Howie}.  A function $t\in \mathcal{T}_S$ is
\emph{well-defined} in $S$  if for any $w, w'\in S$ such that $w =
w'$, one has $t(w) = t(w')$.

 The set of well-defined functions from $\mathcal{T}_S$ can be
utilized to deliver \emph{diffusion} in $S$, i.e., to dissemble an
element of the platform $S$ before transmission by using its
defining relations. If a subset $T\subseteq \mathcal{T}_S$
consisting of well-defined functions acts on $S$, say,
\[
f:S \times T \longrightarrow S \quad \textrm{described by} \quad
f: (w,t)\longmapsto t(w),
\]
such that  recovering  $w$ from $t(w)= f(w,t)$ is computationally
infeasible, then the action $f$ satisfies the principal
requirement of a one-way function.

 A particular example of such a subset $T\subseteq\mathcal{T}_S$
would be $End\,S$, the set (which is actually a monoid)  of
endomorphisms of $S$. Let $G$ be an arbitrary semigroup and let
$\rho: G \rightarrow End\,S$ be a morphism. Then $\rho$ determines
an action of $S$ by its image (denoted by $Im\,\rho$), i.e., $g
\mapsto (t \mapsto t(w))$, for $g\in G$ and $t\in End\,S$.  The
function
\[
f: S \times T \longrightarrow S \quad \textrm{given by}  \quad\
(w,t) \longmapsto t(w) = w'
\]
explicitly defines the action, where $T = Im\,\rho$. If the search
for a $t\in T$ such that $f(w,t) = w'$ is computationally
infeasible, then the action $f$ is an intrinsic one-way function
inherited by $S$ via a semigroup $T$.

 Therefore, an algebraic characterization of a one-way function can
be determined through an action, as specified above, by 
algebraic properties of $S$ and $T$. Without loss of generality,
given  feasibly computable algebraic structures $S$ and $T$, if
there exists an action $f: S\times T \rightarrow S$ such that $f$
is a one-way, well-defined function for fixed values of $T$, then
the elements of $S$ can be manipulated for public-key encryption.

\begin{definition}\label{absPKC}
An algebraic public-key cryptographic system is a tuple $(S,\, T
,\, f;\,\mathcal{H},\, h)$, satisfying the following properties:
\begin{itemize}
     
    \item $S$ and $T$ are feasibly computable algebraic structures
    (e.g. semigroups).
     
    \item $f : S \times T  \rightarrow S$ is an action that is
    one-way and well-defined    for fixed values of $\,T$: given a
    private $t\in T$ and any public    $w \in S$, it is infeasible
    to determine $t$ from $f(w,\,t)$, and for any $w'\in S$ such
    that $w' = w$, one has $f(w',t) = f(w,t)$.
     
    \item $\mathcal{H}$ is a set of auxiliary feasibly computable
    algebraic structures defined for specific protocols (i.e.
    key exchange, decryption, etc).
     
    \item $\,h: X\times Y \rightarrow X$ is an auxiliary action
    (defined for specific protocols), where $X$ and $Y$ are  
    one of the algebraic structures $S,\,T$, or $H\in\mathcal{H}$.
    \end{itemize}
\end{definition}

 Let us now assume $S$ to be a feasibly computable group $G$. With the
developed analysis, we interpret the standard theory of PKC as
arising  from a \emph{permutation representation}
$\rho:G\rightarrow \mathcal{S}_G$, where $\mathcal{S}_G\subseteq
\mathcal{T}_G$ is the symmetric group of $G$, such that $\rho: x
\mapsto x^{\,\rho}$ and $Im\,\rho$ is a subgroup of the group of
automorphisms of $G$, denoted by $Aut\,G$. Since the elements
$\,x^{\,\rho} \in Im \,\rho$ are automorphisms, $\,x^{\,\rho}$
acts by permuting words $g\in G$ with  the capability of providing
cryptographic confusion and diffusion (see e.g. \cite{Garrett}).

 If for every $\,x^{\,\rho} \in Im \,\rho$ the recovery of $g$
from $g\,' = x^{\,\rho}(g)$ is infeasible, then the representation
$\rho$ determines  a one-way function; namely, the group action
\[
f: G\times N \longrightarrow G   \quad \textrm{ defined by } \quad
(g,x^{\,\rho}) \longmapsto x^{\,\rho}(g) = g^{\,\prime}\,,
\vspace{-1mm}
\]
where $N$ is a subgroup of $Aut\,G$. Given
$g\,'=f(g,x^{\,\rho})$, it should be noted that it suffices to
``search" for $(x^{\rho})^{-1}$ in $Aut\,G$ to determine $g\in G$;
this establishes an \emph{automorphism search problem} for $G$.

\begin{definition}\label{problems}
 Let  $F(X)$  be the free group with basis $X$ and let $\langle X ;R
\rangle$ be a presentation of $G$. 
\begin{itemize}
     
    \item Given an arbitrary word $g\in
G$, \emph{the word problem} (WP) is the algorithmic problem of deciding
whether or not $g=1$. 

   \item Given a word $g\in G$ such that $g=1$, the
{\it word search problem} (WSP) is the algorithmic problem of 
``searching" for an explicit expression of $g$ as a product
\,$u_1\,r_1^{\,\epsilon_1}u_1^{-1} \cdots
u_t\,r_t^{\,\epsilon_t}u_t^{-1} = g$, where $u_i\in F(X)$, $r_i\in
R$, and $\epsilon_i \in \{\,\pm 1\,\}$. 

   \item Two words  $g,h\in G$ are
\emph{conjugate} if there is an $x\in G$ such that $xgx^{-1} = h$.
The algorithmic problem of deciding whether or not two arbitrary
words  $g,h\in G$ are conjugate is the \emph{the conjugacy
problem} (CP). 

   \item Given two conjugate words $g,h\in G$, {\it the conjugacy
search problem} (CSP) is the algorithmic problem of ``searching" for an
$x\in G$ satisfying $xgx^{-1} = h$.
    \end{itemize}
\end{definition}

\begin{example}\label{Ex Bn}
\emph{The braid group on $n$ strands, denoted by $B_n$,  with
presentation}
\[
B_n = \langle \,\sigma_1,\ldots,\sigma_{n-1}\; ;\;
\sigma_i\sigma_j\sigma_i=\sigma_j\sigma_i\sigma_j \textrm{ for }
|i-j|=1, \; \sigma_i\sigma_j=\sigma_j\sigma_i \textrm{ for
}|i-j|\geq2\,\rangle\,,
\]
\emph{has the word  problem  solvable in quadratic time. The
braid group $B_n$ is a suggested group-theoretic platform for the
implementation of the conjugacy search problem, see  \cite{AAG}, \cite{KLCHKP}}. 

\emph{We note that the group $Aut\,B_n$ is equal to
$\langle \,Inn \,B_n ,\,\eta\,\rangle$, where $Inn \,B_n$ is the
 group of inner automorphisms of $B_n$ and
$\eta:\sigma\mapsto\sigma^{-1}$, for any $\sigma\in B_n$. Thus,
the general automorphism search problem for $B_n$ basically
reduces to the inner-automorphism search problem for $B_n$, i.e.,
to the conjugacy search problem  for braid groups. }
\end{example}


\section{Commuting Action Key Exchange (CAKE)}
\label{sec CAKE}

 To change the standard methodology of working implicitly just with
the automorphism group of $G$, we generalize an action $f:G\times
Aut \,G\rightarrow G$, using Definition \ref{absPKC}, to a
well-defined action on an algebraic structure $S$ by an algebraic 
structure $N$ for fixed values of $N$. To manifest the advantage
of the abstraction, we construct an algebraic PKC for the
implementation of a key exchange protocol based on a generalization of 
the discrete logarithm problem:

\begin{definition}
[Commuting Action Key Exchange, CAKE]\label{CAKE} Select the
platforms $S$ and $T$  to establish an algebraic PKC tuple
$(S,\,T,\,f;\,\mathcal{H})$, where the auxiliary set $\mathcal{H}$
is $\{ A, B \subseteq T \,|\, \forall \, \alpha \, \in \, A 
~\forall \beta \, \in \, B  ~\alpha\,\beta= \beta\,\alpha\, 
\}$. The key exchange protocol is set for two entities,   Alice
and Bob.

\medskip \emph{\textbf{Protocol:}}

\begin{enumerate}
    
    \item The semigroup $S$, a word $w \in S$,  and a generating set
    for each semigroup in $\mathcal{H}$ are made public.
   
    \item Alice chooses a private word $\alpha \in A$
    satisfying $f(w,\alpha)\neq 1$ and transmits  
      $f(w,\alpha) = w\alpha$ to Bob.
     
    \item Bob chooses a private word $\,\beta \in B\,$
    satisfying $\,f(w,\beta)\neq 1\,$ and transmits
      $f(w,\beta) = w\beta$ to Alice.
    
    \item Alice computes    $f(w\beta,\alpha) = w \,\beta\alpha$ and
    Bob computes $f(w\alpha,\beta) = w \,\alpha \beta$.
    Both entities  establish
    $\:w\,\alpha\beta = w\,\beta\alpha\:$
    as the common secret key.
\end{enumerate}
\end{definition}

\begin{example}
\emph{The Diffie-Hellman protocol becomes an instance of the CAKE
protocol if  the multiplicative group of integers modulo a prime
number  and its standard automorphism group are the chosen
platforms.}
\end{example}

 A simple, well-studied associative algebraic system $S$ with a
single binary operation and a commutative semigroup $T \subseteq
End\,S$ generated by a large set of elements are good potential
candidates for the implementation of CAKE. In this case, both
$\alpha$ and $\beta$ are endomorphisms of $S$, and $\alpha(\beta
(w))=\beta(\alpha(w))$ becomes the common key. Similarly, one can
also use a commutative subsemigroup $T$ of the full transformation
semigroup $\mathcal{T}_S$ containing well-defined functions
$\alpha, \beta$ (not necessarily endomorphisms) of $S$. A basic
example of that kind was given in \cite{Ko}.

\begin{example}
\emph{Let $A, B \subseteq S$ be two subsemigroups of a semigroup
$S$ such that $ab=ba$ for any $a \in A, ~b \in B$. Given a public
element $w \in S$, Alice computes $w\mapsto a_1wa_2$, where $a_1,
a_2\in A$ are her private elements, and transmits this new element
to Bob (after disguising it somehow). Similarly, Bob transmits
$w\mapsto b_1wb_2$, where $b_1, b_2\in B$ are his private
elements. The common key now is $a_1b_1wb_2a_2=b_1a_1wa_2b_2$. }

\emph{Note that if $A, B \subseteq S$ are groups,  the protocol of
Ko, Lee et. al. \cite{KLCHKP} can also be obtained as a special
case of the above protocol where $a_2=a_1^{-1}$ and
$b_2=b_1^{-1}$. }
\end{example}


\section{Classes of groups vs. particular groups}
\label{classes}

 Let $\mathcal{B}$ be the class of braid groups. A
generic element $B_n$ from this class can be chosen from
$\mathcal{B}$ simply by randomly selecting a natural
number for the variable $n$. For general applications, once a
choice for a braid group $B_n$ is made, an algorithm that applies
to this group also applies to other braid groups. Informally
speaking, the braid groups are ``algorithmically homogeneous" and
this can be a drawback for cryptographic applications, as stated
in the Introduction. In the following sections, we address this
issue by considering wider classes of groups.

 In particular, we introduce additional randomness to an algebraic
PKC protocol, requiring that its platforms  be selected at random
from a wider class of groups at the beginning of the generation of
keys. Moreover, isomorphic groups from a wider class provide a
mechanism for diffusion, as examined in the last section of the
paper. The use of isomorphic groups and random selections from a
class of groups is a familiar scenario for cryptosystems; both in
the RSA and in the discrete logarithm cryptosystems, primes are
randomly selected for application, i.e., a multiplicative group of
integers and a subgroup of its automorphism group are randomly
selected.

To exemplify these ideas, we first consider the class of Artin
groups of extra large type for the implementation of the Commuting
Action Key Exchange protocol, via endomorphisms. Second, despite
our belief that Question 1 in the Introduction is likely to have a
negative answer, we give the conjugacy search problem (CSP) a 
benefit of the doubt; we consider the class of groups satisfying
small cancellation conditions $C(4)$ and $T(4)$, but not
$C'(\frac{1}{6})$ (to try to avoid hyperbolic groups), for the
implementation of a cryptosystem relying on CSP. Furthermore,
these classes of groups offer additional properties that can  be
utilized in other algebraic PKC protocols.


\section{The class of Artin groups of extra large type}
\label{Class of Artin Groups}

Let $G\Gamma$ be a group with presentation
\[
G\Gamma = \langle\,g_1,\ldots , g_n\; ; \;
r(g_i,g_j)=1\;\,(\textsl{for} \;\;1 \leq i , j \leq n\;\textsl{
and }\;i \neq j) \; \rangle\, ,
\]
where $n\geq 2$ and $r(g_i,g_j)=1$ is a relator involving two
generators. Given $G\Gamma$ there is an associated labeled graph
$\Gamma$ and vice versa. The vertices of the graph $\Gamma$ are
labeled by the generators of $G\Gamma$. Any two vertices $g_i,
g_j\in \Gamma$ are connected by an edge if there is a relation
$r(g_i,g_j)\in G\Gamma$ between the corresponding generators; in
other words, edges are labeled by relations.

\begin{example}\label{Artin group}
\emph{An \emph{Artin group} $A\Gamma$  is a group with
presentation}
\[
A\Gamma = \langle\;a_1,\ldots,a_n  \; ;\; \mu_{ij} = \mu_{ji}
\;\; \textrm{for} \;\;1 \leq i < j \leq n)\; \rangle\,, \quad
\emph{where }\; \mu_{ij} =
\underbrace{a_i\,a_j\,a_i\ldots}_{m_{ij}}
\]
\emph{and $m_{ij} = m_{ji}$. Artin groups arise as generalizations
of  braid groups, see e.g. \cite{AS}. For an Artin group
$A\Gamma$, the associated labeled graph $\Gamma$ has no multiple
edges or loops. The vertices $a_i$ of $\Gamma$ are the generators
of the Artin group. Any two vertices $a_i, a_j\in \Gamma$ are
connected by an edge, labeled with the integer $m_{ij}\,$,
associated to the relation $\mu_{ij} = \mu_{ji}$ (between the
corresponding generators $a_i, a_j\in A\Gamma$).}
\end{example}

 In general, automorphisms (or endomorphisms) of the graph $\Gamma$
induce automorphisms (or endomorphisms) of the group $G\Gamma$.
Therefore, the graph associated to $G\Gamma$ gives us a direct
procedure for the construction of a semigroup $T \subseteq End
\,G\Gamma$ that can contain a large pool of commuting elements.
This is a necessary condition for common key extraction by
legitimate parties in the application of the Commuting Action Key
Exchange protocol (CAKE, Definition \ref{CAKE}). To construct the
corresponding semigroup $T$ with sufficiently many endomorphisms,
a graph $\Gamma$ can be chosen to be a tree. The procedure
implemented for the Ko-Lee protocol can then be utilized to
provide for commuting endomorphisms, i.e., one splits the vertices
of the graph into two disjoint sets such that each of the
entities, Alice and Bob, select endomorphisms which act on their
own set.

\begin{example}
\emph{The relations of the braid groups $B_n$ involve two
generators. The corresponding graph associated to $B_n$ is just a
simple path, and it  has only one automorphism that induces the
following automorphism of $B_n$:
 $\sigma_i \mapsto \sigma_{n-i}$, which  happens
 to be an  inner automorphism of $B_n$.  For other $G\Gamma$
groups, however, their corresponding graphs are more complex,
 and it is easy to arrange for a large
semigroup (or a group) $T \subseteq End\,G\Gamma$ of endomorphisms
(or automorphisms).}
\end{example}

  Artin groups $A\Gamma$ with the property that all the integers
$m_{ij} \geq 4$\, are  called \emph{Artin groups of extra large
type}. A  tree $\Gamma$ can be associated to an Artin group of
extra large type, providing a direct procedure for   constructing
a semigroup $T \subseteq End \,A\Gamma$. Moreover, Artin groups of
extra large type are automatic \cite{Peifer}, thus the word
problem for  groups in  this class can be solved in quadratic
time, and by a result of \cite{KMSS2}, the word problem is
solvable in linear time on average. Therefore, we can  suggest the
class of Artin groups of extra large type  as platforms for CAKE.


\subsection{Key exchange protocol based on Artin groups}

 In this section we present the class of Artin groups of extra
large type as an implementable class for CAKE.

\smallskip

\noindent\textbf{Key generation:} Randomly select a finite rooted
tree $\Gamma$ with $l$  levels such that the degree of the root is
equal to 2, and the degrees of all other vertices are between 2
and an integer $m$,  with the exception of the end vertices whose
degrees are 1. Associate to the tree $\Gamma$ an Artin group
$A\Gamma$ of extra large type by labelling each vertex of $\Gamma$
with a letter $a_i$ and numbering an edge by a (random)
$m_{ij}\geq 4$ if there are  two corresponding vertices $a_i$ and
$a_j$ incident to this edge.

 Let $a_k$ be the root of the tree and let $\Gamma_0=
\Gamma - a_k$ be the   subgraph obtained by deleting
the root $a_k$. The graph $\Gamma_0$ consists of two finite
disjoint subtrees, say, $\Gamma_{A}$ and $\Gamma_{B}$, that are
spliced by the root $a_k$. The associated subgroups are $A\Gamma_{A}$
and $A\Gamma_{B}$.

 The sets of graph endomorphisms of $\Gamma_{A}$ and $\Gamma_{B}$
induce the submonoid of endomorphisms $End\,A\Gamma_{A} \times
End\,A\Gamma_{B}\subseteq End\,A\Gamma$ such that for any
$\alpha\in End\,A\Gamma_{A}$ and $\beta\in End\,A\Gamma_{B}$ both
$\alpha$ and $\beta$ commute: $\alpha\,\beta = \beta\,\alpha$. In
order for both submonoids to act non-trivially on a public word
$w\in A\Gamma$, the word must involve  some generating elements
$a_1,\ldots,a_p \in A\Gamma_{A}$ and   some generating elements
$b_1,\ldots,b_q \in A\Gamma_{B}$, i.e.,  $w=w(a_1,\ldots,a_p,\,b
_1,\ldots,b_q)$.

\medskip

\noindent {\large\textbf{CAKE for Artin groups of extra
large type.}} Choose a random  Artin group $A\Gamma$  of extra large type
to be the   platform $S$ for the CAKE tuple
$(S,\,T,\,f;\,\mathcal{H})$,  and let $T= End\,A\Gamma_{A} \times
End\,A\Gamma_{B}$. Define $\mathcal{H}$ to be the set
$\{End\,A\Gamma_{A} ,\, End\,A\Gamma_{B}\}$. The protocol is set
for Alice and Bob.

\medskip

\textbf{Protocol:}

\begin{enumerate}
     
    \item The random group $A\Gamma$, a word
    $w=w(a_1,\ldots,a_p,\,b _1,\ldots,b_q)\in A\Gamma$ and a generating
    set for each element of $\mathcal{H}$ are made public.
     
    \item Alice chooses a private word $\alpha \in
    End\,A\Gamma_{A}$ and transmits $f(w,\alpha) = w^{\,\alpha}$ to Bob.
     
    \item Bob chooses a  private word $\beta \in
    End\,A\Gamma_{B}$ and transmits $f(w,\beta) = w^{\,\beta}$ to Alice.
     
    \item Alice computes $f( w^{\,\beta},\alpha) = w^{\,\beta\,\alpha}$
    and Bob computes $f(w^{\,\alpha},\beta) = w^{\,\alpha\,\beta}$.
    Alice and Bob set
    \[
    \; w^{\,\alpha\,\beta}= w^{\,\beta\,\alpha}
    \]
    as their common secret key.
\end{enumerate}

\noindent \textbf{Remark.} By introducing randomness in the
selection of the group $A\Gamma$, we make the present approach
dynamic. The class of  Artin groups of extra large type seems to
be less ``algorithmically homogeneous" than, say,  the class of
braid groups. In general, algorithmic non-homogeneity can disrupt
  general algorithmic methods an opponent might obtain for the
purpose of acquiring a private key. For example, a typical
endomorphism (non-automorphism) for $A\Gamma$ would be merging two
terminal children vertices of the same parent, ``confusing" the
length of the word $w$. As a result, the effect of such an
endomorphism  on the length of a generic element of the group is
no longer predictable, placing  length attacks in question.


\section{A class of small cancellation groups}
\label{small cancellation}

 In this section, we follow Lyndon and Schupp \cite{L-S}. For facts
about small cancellation theory the reader is referred to this
source for further reading. Let $F(X)$ be the free group with a
basis $X = \{\, x_i \,|\, i\in I \,\}$, where $I$ is an indexing
set. Let $\epsilon_k\in \{\pm 1\}$, where $1\leq k\leq n$. A word
$w(x_1,\ldots,x_n)=x_{i_1}^{\epsilon_1}x_{i_2}^{\epsilon_2}\cdots
x_{i_n}^{\epsilon_n}$ in $F(X)$, with all   $x_{i_k}$ not
necessarily distinct, is a \emph{reduced $X$-word} if
$x_{i_k}^{\epsilon_k}\neq x_{i_{k+1}}^{-\,\epsilon_{k+1}}$, for
$1\leq k \leq n-1$. In addition, the word $w(x_1,\ldots,x_n)$ is
\emph{cyclically reduced} if it is a reduced $X$-word and
$x_{i_1}^{\epsilon_1} \neq x_{i_{n}}^{-\,\epsilon_n}$. A set $R$
containing cyclically reduced words from $F(X)$ is
\emph{symmetrized} if it is closed under cyclic permutations and
taking inverses.

 Let $G$ be a group  with presentation $\langle X;R\rangle$. A
non-empty word $u \in F(X)$ is called a \emph{piece} if there are
two distinct relators $r_1, r_2 \in R$ of $G$ such that $r_1 = u
v_1$ and $r_2 = u v_2$. The group $G$  belongs to the class  
  $C(p)$   if no element of
$R$ is a product of fewer than $p$ pieces. Also, the group $G$
belongs to the class  
$C'(\lambda)$   if  for every $r\in R$ such that $r = uv$ and
$u$ is a piece, one has  $|u| < \lambda |r|$.

 In particular, if $G$ belongs to the class    $C'(\frac{1}{6})$, 
then Dehn's algorithm solves the word problem  for $G$. Thus, if $G$ is a  
finitely presented group from the class  $C'(\frac{1}{6})$, then it is
\emph{hyperbolic}.

 \begin{example}\label{Ex c4 not 1-6}
\emph{ Let $\langle\,x_1,x_2,x_3 \,;\, x_1^2x_2x_3^2x_2^{-1}=1, \,
x_2^2x_3x_1^2x_3^{-1}=1\,\rangle \,$ be a presentation of a group
$G$. Now, $x_1^{\pm 2},\,x_2^{\pm 1},\,x_2^{\pm 2},\,x_3^{\pm
1},\,x_3^{\pm 2},\,(x_2x_3)^{\pm 1}$ and $(x_2x_3^{-1})^{\pm 1}$
are the pieces of $G$, and every relator is a 
  product of four of these pieces. Therefore, the group $G$ is
in the class of $C(4)$ groups. However, $G$ 
is not in the class of $C'(\frac{1}{6})$; for $i=1,\,2$ and $3$,
the pieces $x_i^{\pm 1}$, satisfy the property $|x_i^{\pm 1}| =
\frac{1}{6} |x_1x_2x_3x_4^2x_2^{-1}|$ and $|x_i^{\pm 1}| =
\frac{1}{6} |x_2^2x_3x_1x_4x_3^{-1}|$. }
\end{example}

 The solution of the conjugacy problem (CP) is  irrelevant
for the implementation of a cryptographic protocol utilizing the
computational difficulty of the conjugacy search problem (CSP).
However,  reasonable evidence of a potentially computationally
hard CSP is  provided  if there is no known polynomial time
algorithm for CP.

 For a class of small cancellation groups possessing  the property
of ``no known polynomial time algorithm for CP", we need one more
condition. A group $G$ with finite presentation $\langle
X;R\rangle$ belongs to the class $T(q)$ for a natural number $q$
 if  for any
sequence $r_1,\ldots,r_n \in R$, with $3\leq n < q$ and $r_i\neq
r_{i+1}^{-1}$,   at least one of the products $r_1
r_2,\ldots,\,r_{n-1}r_n,\,r_n r_1$ is cyclically reduced without
cancellation.

A group $G$ with presentation $\langle\,X\,|\,R\,\rangle$ is said
to be a small cancellation group of type $C(p)$-$T(q)$ if  it   
  belongs to the classes $C(p)$ and   $T(q)$.
By Theorem V.6.3 of \cite{L-S}, the word problem is
solvable in the class of small cancellation groups of type
$C(4)$-$T(4)$. If hyperbolic groups $C'(\frac{1}{6})$ are avoided, 
then, generally, there is no known polynomial
time algorithm for solving the conjugacy search problem for groups in 
  this class (even though the conjugacy problem is solvable by 
\cite[Theorem V.7.6]{L-S}). 

 Thus, in this class,
legitimate entities can choose a random group and implement an
algebraic PKC  protocol, e.g. CAKE, that relies on the hardness
of the conjugacy search problem or a harder problem that can
potentially arise (as we indicated in Sections \ref{sec Alg PKC}
and \ref{sec CAKE}). 

\begin{example}\label{c4-t4 not 1-6}
\emph{Consider the presentation $\langle\,x_1,x_2,x_3 \,;\,
x_1^2x_2x_3^2x_2^{-1}, \, x_2^2x_3x_1^2x_3^{-1}\,\rangle\,$ for
$G$, the group of Example \ref{Ex c4 not 1-6}. For any $r_1,\,r_2$
and $r_3$, no two of which are inverse of one another, from the
symmetrized set $\{x_1^2x_2x_3^2x_2^{-1}, \,
x_2^2x_3x_1^2x_3^{-1}\}$, no cancellation is possible in at least
one of the words $r_1 r_2,\,r_2r_3$ and $r_3r_1$. Therefore, $G$
belongs to the class of $T(4)$ and  $C(4)$ groups, but not the
class of $C'(\frac{1}{6})$ groups. }
\end{example}


\section{Diffusion}
\label{diffusion}

 In this section, we offer a method that can, in our opinion,
substantially enhance the ``diffusion", i.e., the process of
disguising an element of a given group by using defining
relations. This method is not brand new, but it was used before in
a different context, namely, in attempts to attack the
Andrews-Curtis conjecture, a notoriously difficult problem in
low-dimensional topology and combinatorial group theory (see e.g.
\cite{MMS}).

 The idea is to break down defining relations of a
group into ``small pieces". More formally, we replace a given
group $G$ by an isomorphic group where all relators have length at
most 3. Intuitively, diffusion should be  easier to achieve in
groups with  shorter defining relations, so we hope that our idea
is useful.

 The procedure itself is quite simple. Let $G$ have a presentation
$\langle x_1,...,x_n; r_1,...,r_k\rangle$ in terms of generators
$x_1,...,x_n$  and defining relations $r_1,...,r_k$. We are going
to obtain a different presentation for $G$ by using {\it Tietze
transformations}
 (see e.g. \cite{L-S}); these are elementary isomorphism-preserving
operations on presentations of groups.

 Specifically, let, say,  $r_1=x_i x_j u$, $1 \le i,j \le n$.
We introduce a new generator $x_{n+1}$  and a new relator
$r_{k+1}=x_{n+1}^{-1}x_i x_j$. The group with the presentation
$\langle x_1,...,x_n, x_{n+1}; r_1,...,r_k, r_{k+1}\rangle$ is
obviously isomorphic to $G$. Now if we replace $r_1$ with
$r_1'=x_{n+1} u$, then the presentation $\langle x_1,...,x_n,
x_{n+1}; r_1',...,r_k, r_{k+1}\rangle$ will again define a group
isomorphic to $G$, but now the length of one of the defining
relations ($r_1$) has decreased by 1. Continuing in this manner,
we can eventually obtain a presentation where all relators have
length at most 3, at the expense of introducing more generators.

 Apparently, relators of length at most 3 can provide a very good diffusion,
but the natural question now is: why cannot the opponent convert
the new presentation back to the original one  and take it from
there? This, indeed, may work with some of the protocols, but let
us have a look at the situation where applying an endomorphism of
a group to an element is involved.

 Suppose a group $G'$ is isomorphic to a  group $G$ in the way described above.
 Let $w' \in G'$, and let
$\varphi$ be an endomorphism of $G'$ applied to $w'$. The opponent
can convert $w'$  and $\varphi(w')$ to elements $w$  and $u$,
respectively, of the group $G$, by using relations of the form
$x_s=x_{i_1}^{\pm 1} x_{i_2}^{\pm 1}$, where $x_s$ are ``new" generators
and $x_{i_1}, x_{i_2}$ are  ``old" generators. Then the opponent
may try to find an endomorphism $\psi$  of $G$ such that
$u=\psi(w)$ as follows.

 Suppose we know that
$\varphi$ takes generators $x_i'$  of the group $G'$ to some
$y_i'$. An obvious way to ``lift" $\varphi$    to an endomorphism
of $G$ would be to convert $y_i'$ to $y_i \in G$ (again,  by using
relations of the form $x_s=x_{i_1}^{\pm 1} x_{i_2}^{\pm 1}$), then let
$\psi$ be the mapping of  $G$ that takes $x_i$ to $y_i$.

 This however may not work (and typically will not work) because the
endomorphism $\varphi$, restricted to the  ``old" generators
(i.e., to the generators of $G$) may not respect the original
relations of the group $G$. We can therefore have an element $u
\in G$ such that $\psi(w)=u$ in the group $G$, but $\varphi(w')\ne
u$ in the group $G'$.

 The only way to properly ``lift" $\varphi$ to an endomorphism of $G$ would
be to combine it with an isomorphism $f: G' \to G$, but the latter
is by no means easy to explicitly compute, even if the whole chain
of Tietze transformations is known to the opponent, which  does
not have to be the case. Incidentally, neither has the original
group $G$ to be known to the public.


\medskip
\noindent
 Department of Mathematics, The City  College  of New York, New York,
NY 10031

\medskip

\noindent {\it e-mail addresses\/}:
shpil@groups.sci.ccny.cuny.edu,
  nyzapata@verizon.net
\medskip

\noindent {\it http://www.sci.ccny.cuny.edu/\~\/shpil/   }

\end{document}